\documentclass[reqno,11pt]{amsart}
\usepackage{amssymb}
\usepackage{mathrsfs}
\usepackage{cite}
\usepackage{amsfonts}
\usepackage{amsthm}
\usepackage{amsgen}
\usepackage{amsmath}

\def\wh{\widehat}

\def\pv#1{\ensuremath{{\bf#1}}}

\def\ilim{\varprojlim}

\def\J{\mathrel{{\mathscr J}}} 

\def\H{\mathrel{{\mathscr H}}} 

\def\e<{\leq _{E}}

\def\NN{\ensuremath{\mathbb N}}
\def\malce{\mathbin{\hbox{$\bigcirc$\rlap{\kern-8.3pt\raise0,50pt\hbox{$\mathtt{m}$}}}}}

\def\1sk{^{(1)}}

\def\to{\rightarrow}

\def\Thmname{Theorem}
\def\Propname{Proposition}
\def\Lemmaname{Lemma}
\def\Definitionname{Definition}
\newtheorem{Thm}{\Thmname}

\newtheorem{Lemma}[Thm]{\Lemmaname}
{\theoremstyle{definition}
}
{\theoremstyle{remark}
}
\newtheorem{Cor}[Thm]{Corollary}

{\theoremstyle{remark}
}

\numberwithin{equation}{section}

\title{A Combinatorial Property of Ideals in Free Profinite Monoids}

\author{Benjamin Steinberg}
\address{School of Mathematics and
  Statistics\\ Carleton
University\\ Ottawa, Ontario K1S 5B6\\ Canada}
\email{bsteinbg@math.carleton.ca}
\thanks{The author gratefully acknowledge the support of NSERC}
\date{November 8, 2008}

\begin{document}
\maketitle

The reader is referred to~\cite{qtheor} for all undefined notation concerning finite and profinite semigroups.
We assume throughout this note that $\pv V$ is a pseudovariety of monoids~\cite{qtheor,Almeida:book,Eilenberg}
closed under Mal'cev product with the pseudovariety $\pv A$ of
aperiodic monoids, i.e., $\pv A\malce \pv V=\pv V$.  Denote by
$\wh{F}_{\pv V}(A)$ the free pro-$\pv V$ monoid on a profinite space
$A$~\cite{Almeida:book,qtheor,AWsurvey}.  In this note we prove the following theorem.

\begin{Thm}\label{main}
Suppose that $\alpha_1,\ldots,\alpha_m\in \wh{F}_{\pv V}(A)$ and
$I_1,\ldots,I_n$ are closed ideals in $\wh{F}_{\pv V}(A)$ where $m\leq n$. If
$\alpha_1\cdots\alpha_m\in I_1\cdots I_n$, then $\alpha_i\in I_j$ for
some $i$ and $j$.
\end{Thm}

Before proving the  theorem, we state a number of consequences.  Recall
that an ideal $I$ in a semigroup is \emph{prime} if $ab\in I$
implies $a\in I$ or $b\in I$.

\begin{Cor}\label{primeideal}
Let $I=I^2$ be a closed idempotent ideal of $\wh{F}_{\pv V}(A)$.  Then $I$
is prime.
\end{Cor}
\begin{proof}
Suppose that $ab\in I=I^2$.  Then $a\in I$ or $b\in I$ by Theorem~\ref{main}.
\end{proof}

An element $a$ of a semigroup $S$ is said to be regular if there
exists $b\in S$ so that $aba=a$.  Any regular element of a profinite
semigroup generates a closed idempotent ideal.  Hence we have:

\begin{Cor}\label{MishaJorge}
Every regular element of $\wh{F}_{\pv V}(A)$ generates a prime ideal.
In particular, the minimal ideal of $\wh{F}_{\pv V}(A)$ is prime.
\end{Cor}

The second statement of Corollary~\ref{MishaJorge} was first proved by
Almeida and Volkov using techniques coming from symbolic
dynamics~\cite{AlmeidaVolkov2}.

Our next result generalizes a result of Rhodes and the author showing
that all elements of finite order in $\wh{F}_{\pv V}(A)$ are group
elements, which played a key role in proving that
such elements are in fact idempotent~\cite{projective}.
Let $\wh{\NN}$ denote the profinite
completion of the monoid of natural numbers; it is in fact a profinite semring.  We use $\omega$ for the
non-zero idempotent of $\wh{\NN}$.

\begin{Cor}
Let $\alpha\in \wh{F}_{\pv V}(A)$ satisfy $a^n=a^{n+\lambda}$ for some
positive integer $n$ and some $0\neq \lambda\in \wh{\NN}$.  Then $a$ is a
group element, i.e., $a=a^{\omega}a$.
\end{Cor}
\begin{proof}
It is immediate that $a^{n} = a^{n+k\lambda}$ for all $k\in \mathbb N$ and hence all $k\in\wh{\mathbb N}$.  Thus $a^n = a^{n+\omega\lambda} = a^{\omega}$.  Since $a^{\omega}$ is regular, it generates a prime ideal by  Corollary~\ref{primeideal}.  It follows that $a\J a^{\omega}$ and hence $a\H
a^{\omega}$. Thus $a$ is a group element, as required.
\end{proof}

To prove the theorem, we use the Henckell-Sch\"utzenberger expansion.
Let $M$ be a finite monoid generated by a set $A$.  For any element
$\alpha\in \wh{F}_{\pv V}(A)$, we write $[\alpha]_M$ for its image in
$M$.   For $w\in A^*$ (the free
monoid on $A$), define
$\mathrm{cut}_n(w)$ to be the set of all $n$-tuples $(m_1,\ldots,m_n)$
of $M$ such that there exists a factorization $w=w_1\cdots w_n$ with
$[w_i]_M=m_i$, for $i=1,\ldots,n$.  It is well known that the
equivalence relation on $A^*$ given by $u\sim v$ if
$\mathrm{cut}_n(u)=\mathrm{cut}_n(v)$ is a congruence of finite index
contained in the kernel of the natural map $A^*\to
M$~\cite{BR--exp,ourstablepairs}.  Moreover, if we denote by $M^{(n)}$
the quotient $A^*/{\sim}$, then the natural map $\eta\colon M^{(n)}\to
M$ is aperiodic and hence if $M\in \pv V$, then $M^{(n)}\in \pv
V$ under our hypothesis on $\pv V$~\cite{BR--exp,ourstablepairs}.

Our main theorem relies on the following simple factorization lemma for free monoids.

\begin{Lemma}\label{factor}
Let $w\in A^*$ and suppose $w=u_1\cdots u_m=v_1\cdots
v_n$ where $m\leq n$.  Then $v_j$ is a factor of $u_i$ for
some $1\leq i\leq m$ and $1\leq j\leq n$.
\end{Lemma}
\begin{proof}
If any $v_j$ is empty, we are done so assume now each $v_j$ is non-empty and
hence $w$ is non-empty. If any $u_i$ is empty we
may omit it, so assume $u_i$ is non-empty for $i=1,\ldots,m$.
Define a function $f\colon \{1,\ldots,n\}\to \{1,\ldots,m\}$ by
$f(j)=i$ if the last letter of $v_j$ belongs to $u_i$.  The map $f$ is
monotone.  Suppose first that $f$ is not injective.  Then there exists
$2\leq j\leq n$ so that $f(j-1)=f(j)$.  In this case, $v_j$ is a
factor of $u_i$ where $i=f(j)$.  If $f$ is injective, then since it is
monotone and
$m\leq n$, we must have that $m=n$ and $f$ is the identity map.  In
this case, $u_1$ is a factor of $v_1$.
\end{proof}

\begin{proof}[Proof of Theorem~\ref{main}]
Suppose first that $A$ is finite and  that
$\alpha_2,\ldots,\alpha_m\notin I_j$ for any $1\leq j\leq n$ and that
$\alpha_1\notin I_2,\ldots, I_n$.  We show that $\alpha_1\in I_1$.
Since $I_1$ is closed, it suffices to show that $\pi(\alpha_1)\in
\pi(I_1)$ for all continuous surjective homomorphisms $\pi\colon
\wh{F}_{\pv V}(A)\to V$ with $V\in \pv V$.   Using that
$I_1,\ldots,I_n$ are closed we may assume without loss of generality $\pi(\alpha_i)\notin I_j$ for $2\leq i\leq m$ and $1\leq j\leq n$, or $i=1$ and $2\leq j\leq n$.  By assumption  $\alpha_1\cdots\alpha_m\in I_1\cdots
I_n$ so there exist $\beta_j\in I_j$, for $j=1,\ldots,n$, so that
$\alpha_1\cdots\alpha_m=\beta_1\cdots\beta_n$.

Since $A^*$ is dense in $\wh{F}_{\pv V}(A)$, we can find words
$u_1,\ldots,u_m\in A^*$ so that $[u_i]_{V^{(n)}}=[\alpha_i]_{V^{(n)}}$, for
$i=1,\ldots,m$ and words
$w_1,\ldots,w_n$ so that $[w_j]_{V^{(n)}}=[\beta_j]_{V^{(n)}}$ for
$j=1,\ldots,n$.  Then \[[u_1\cdots u_m]_{V^{(n)}}=[\alpha_1\cdots
\alpha_m]_{V^{(n)}}=[\beta_1\cdots \beta_n]_{V^{(n)}}= [w_1\cdots
w_n]_{V^{(n)}}\] and so there exists a factorization $u_1\cdots
u_m=v_1\cdots v_n$ so that $[v_j]_{V}=[w_j]_{V}$, for $i=1,\ldots,n$.
By Lemma~\ref{factor} it follows that there exist $i$ and $j$ so that
$v_j$ is a factor of $u_i$.  Since $[v_j]_V=[w_j]_V=[\beta_j]_V\in
\pi(I_j)$, it follows that $\pi(\alpha_i)=[u_i]_V\in \pi(I_j)$.  By
our assumption on $\pi$, it follows that $i=1=j$. Thus $\pi(\alpha_1)\in
\pi(I_1)$, as required.  This completes the proof when $A$ is finite.

Suppose next that $A$ is profinite; so $A=\ilim_{d\in D} A_d$ with
$D$ a directed set and $A_d$ finite for $d\in D$.
Then (cf.~\cite{Almeida:book,AWsurvey}) one has $\wh{F}_{\pv V}(A) =
\ilim_{d\in D} \wh{F}_{\pv V}(A_d)$.  Suppose now that  $\alpha_2,\ldots,\alpha_m\notin I_j$ for any $1\leq j\leq n$ and
$\alpha_1\notin I_2,\ldots, I_n$.  We show that $\alpha_1\in I_1$.
Since $I_1$ is closed, it suffices to prove, for all $d\in D$, that $\pi_d(\alpha_1)\in
\pi_d(I_1)$ where $\pi_d\colon \wh{F}_{\pv V}(A)\to \wh{F}_{\pv
  V}(A_d)$ is the canonical projection.

Now $\pi_d(I_j)$ is a closed
ideal for $1\leq j\leq n$, so we may assume without loss of generality
that  $\pi_d(\alpha_i)\notin I_j$ for $2\leq i\leq m$ and $1\leq j\leq
n$, or $i=1$ and $2\leq j\leq n$.  By the previous case,
$\pi_d(\alpha_i)\in \pi_d(I_j)$ for some $i$ and $j$.  By assumption
we must have $i=1$ and $j=1$.  This completes the proof.
\end{proof}

\bibliographystyle{abbrv}
\bibliography{standard2}

\end{document}